\documentclass[11pt,a4paper]{article}
\usepackage{a4wide}
\usepackage{latexsym} 
\usepackage{amsmath}
\usepackage{amssymb}
\usepackage{stackrel}  
\usepackage{color}
\usepackage{graphicx}
\usepackage[linesnumbered,ruled,vlined]{algorithm2e}

\usepackage{amsmath, amsthm, amsfonts}
\usepackage{authblk}

\allowdisplaybreaks[3] 

\usepackage{hyperref}

\usepackage[utf8]{inputenc}
\usepackage{color}

\allowdisplaybreaks[3]

\theoremstyle{plain}
\newtheorem{theorem}{Theorem}[section]
\newtheorem{proposition}[theorem]{Proposition}
\newtheorem{lemma}[theorem]{Lemma}

\newtheorem{defin}{Definition}[section]

\theoremstyle{remark}
\newtheorem{remark}[theorem]{Remark}

\numberwithin{equation}{section} 

\newcommand{\C}{\mathbb{C}}

\newcommand{\R}{\mathbb{R}}

\newcommand{\Ea}{E_{\alpha}}
\newcommand{\Eah}{E_{\alpha,h}}
\newcommand{\La}{\Lambda_{\alpha}}

\newcommand{\I}{\mathcal{I}}
\newcommand{\cc}{\gamma}

\begin{document}
\title{Dunkl derivative from moment differentiation}

\author[1]{Edmundo J. Huertas}
\author[2]{Alberto Lastra}
\author[3]{Judit M\'inguez Ceniceros}
\affil[1]{Universidad de Alcal\'a, Dpto. F\'isica y Matem\'aticas, Alcal\'a de Henares, Madrid, Spain. {\tt edmundo.huertas@uah.es}}
\affil[2]{Universidad de Alcal\'a, Dpto. F\'isica y Matem\'aticas, Alcal\'a de Henares, Madrid, Spain. {\tt alberto.lastra@uah.es}}
\affil[3]{Universidad de La Rioja, Dpto. Matem\'aticas y Computaci\'on, Logro\~no, Spain. {\tt judit.minguez@unirioja.es}}

\maketitle
\thispagestyle{empty}
{ \small \begin{center}
{\bf Abstract}
\end{center}

The work analyzes the theory of Dunkl operator as a moment differential operator. This last operator generalizes the first one whenever the sequence of moments satisfies appropriate classical properties, classically considered in the general theory of ultraholomorphic and ultradifferentiable classes of functions. In this sense, the theory of Dunkl operator is then generalized. On the other hand, some features developed in Dunkl theory, such as Dunkl translation, have not been considered in the theory of moment differential equations yet, which leads to a common mutualism involving both theories.

\smallskip

}
\bigskip \bigskip

\section{Introduction}

The main aim in the present work is to apply the moment differentiation and the theory of summability to the Dunkl operator introduced in~\cite{Du}. Up to our knowledge, both theories have been developed independently. However, each of them can take advantage of the other.

The purpose of this seminal research is to describe the connection points of both theories, and provide a point of departure for a joint future research.

On the one hand, the theory of moment differential equations is a quite recent theory whose origin can be found in the work~\cite{bayo}, by W. Balser and M. Yoshino, where the authors put forward the concept of moment differentiation $\partial_m$, departing from a fixed sequence of positive real numbers $m$. In principle, moment differentiation is defined as a formal operator acting on formal power series and can be naturally extended to holomorphic functions defined on some neighborhood of the origin. In recent years, moment differentiation has been extended to holomorphic functions defined on sectors of the complex domain with vertex at the origin determining the generalized sum (see~\cite{lamisu}) or the generalized multisum  (see~\cite{lamisu4}) of a formal power series. The development of the theory on the analytic and formal solutions to functional equations involving moment derivation has been notorious during the last decade such as in~\cite{lamisu2,lamisu3,mi,michalik12,misutk,mitk}. We also refer to the recent works~\cite{laspri,lastra} regarding the solutions to linear systems of moment differential equations, the works on moment differential operator preserving summability properties~\cite{ichmich,lasmich,misutk}, and the book~\cite{remybook}, and the references therein. Due to the desirable properties of the ultraholomorphic and ultradifferentiable classes of functions related to some sequence $m$ as above, the sequence is usually considered to be a strongly regular sequence. A close branch of research deals with the properties of the previous spaces of functions and the formal representation of their elements near the origin. In this direction, we refer to the recent advances made in~\cite{jimisasc, jisasc}, and other references mentioned in the sequel.

The main advantage of the theory of moment differential equations is its versatility. A modification of the sequence $m$ leads to different concrete realizations of the moment derivative to a usual derivative, Caputo fractional derivative, $q$-derivative, $(p,q)$-derivative, etc... Dunkl operator shows up as another realization of moment derivative when fixing the moment sequence to be the sequence of Dunkl factorials. Given a sequence of positive real numbers $m$, and considering the moment derivative $\partial_m$, a linear system of moment differential equations can be stated:
\begin{equation}\label{eq:moment-syst}
\partial_my(z)=Ay(z),
\end{equation}
where $A\in\C^{n\times n}$ is a constant matrix, $y(z)=(y_1(z),\ldots,y_n(z))^T$ is a vector of unknown functions for some positive integer $n\ge 1$ and
\[
\partial_my(z)=(\partial_m y_1(z),\ldots,\partial_m y_n(z))^T.
\]
Under certain conditions over the sequence $m$, A. Lastra in~\cite{lastra} has obtained the explicit description of the general solutions to a linear system like~\eqref{eq:moment-syst} and the asymptotics of the solutions.

As a matter of fact, the main advances made in the last years regarding functional equations involving moment derivatives, known as moment differential equations, can be applied in this concrete setting of Dunkl operator due to the fact that the sequence of Dunkl factorials turns out to be a strongly regular sequence. In Section~\ref{sec3}, we recall the main advances made on the theory of solutions to linear systems of moment differential equations rephrasing them in Dunkl context. 

On the other hand, being Dunkl operator a particular case of moment differentiation, the theory embracing such operator has far more development, since 1989, appearing in the work~\cite{Du}, by C. Dunkl. We take advantage of this to provide further advances in the theory of moment differential equations, as it is shown in Section~\ref{sec4}, where we study moment translation equations. Moreover, in the last Section we analize the corresponding $m$-even translation operator, defined as the symmetrization of the generalized $m$-translation. This construction isolates the even component of generalized shifts and, consequently, inherits all analytic properties established for the $m$-translation without additional hypotheses, and we also describe its spectral action on generalized $m$-exponentials. This viewpoint helps to clarify the structure of even symmetries in the theory and supplies a tool for analyzing functional equations built from generalized translations.

Therefore, the main goal of this paper is to relate both theories. On the one hand, we will study the solutions and their asymptotics of a linear system of moment differential equations when our moment sequence is the sequence of Dunkl factorials, \eqref{eq:ccna}, obtaining a nice example of a strongly regular sequence. On the other hand, the $m$-translation~\eqref{eq:m-tras} will be defined generalizing Dunkl context, being applied to solve certain families of functional equations.

\section{Dunkl derivative as a moment derivative}

In this section, we recall the definition and main properties of Dunkl derivative and regard it as a moment derivation.

Let $m=\{m(p)\}_{p=0}^\infty$ be a sequence of positive real numbers. Let $\hat{f}(z)=\sum_{p\ge 0}a_p\frac{z^p}{m(p)}$ be a formal power series with complex coefficients. The moment derivative of $\hat{f}(z)$ is given by
$$
\partial_m \hat{f}(z)=\sum_{p\ge 1}a_p\frac{z^{p-1}}{m(p-1)}.
$$ 
Note that when the sequence $m$ is given by $\{p!\}_{p=0}^{\infty}$, the moment derivative of a formal series $\hat{f}(z)=\sum_{p\ge 0}a_p\frac{z^p}{p!}$ is the classical derivative operator
\[
\partial_{p!}\hat{f}=\frac{d}{dx}\hat{f}(x)=\sum_{p\ge 1}a_p\frac{z^{p-1}}{(p-1)!}.
\]
Further realizations of the moment derivative or quite related to it have appeared in the literature, such as Caputo fractional derivative, $q$-derivative, $(p,q)$-derivative, etc. We refer to the detailed examples in~\cite{su}.

Let us now consider the Dunkl factorials. Let $\alpha>-1$ and let  $\gamma=\{\cc_{p,\alpha}\}_{p\ge 0}$ be the sequence of Dunkl factorials defined by
\begin{equation}
\label{eq:ccna}
\cc_{p,\alpha} =
\begin{cases} 
  2^{2k}k!\,(\alpha+1)_k, & \text{if $p=2k$},\\
  2^{2k+1}k!\,(\alpha+1)_{k+1}, & \text{if $p=2k+1$},
\end{cases}
\end{equation}
where $(a)_k$ denotes the Pochhammer symbol
\[
  (a)_k = a(a+1)(a+2) \cdots (a+k-1) = \frac{\Gamma(a+k)}{\Gamma(a)}
\]
(with $k$ a non-negative integer).
In this case, the moment derivative matches with the well-known Dunkl operator (see~\cite{Du, Ros})
$$\partial_{\gamma}f(z)=\La f(z) := \frac{d}{dz}f(z)+\frac{2\alpha+1}{2}\frac{f(z)-f(-z)}{z},$$
 If we denote
\[
\theta_{p,\alpha}:=\frac{\cc_{p,\alpha}}{\cc_{p-1,\alpha}}=p+\frac{2\alpha+1}{2}(1-(-1)^p),
\]
it is easy to prove that
\[
\La x^p=\theta_{p,\alpha}x^{p-1},\quad p\ge 1.
\]
It is worth remarking that the sequence $\{\theta_{p,\alpha}\}_{p=0}^{\infty}$ is also of great importance in the framework of ultraholomorphic and ultradifferentiable classes of functions, known in this context as the sequence of quotients. Many features associated to the previous spaces of functions are determined in terms of properties satisfying the sequence of quotients. See~\cite{sanzrev} and the references therein for more details.

Note that when $\alpha=-1/2$ the Dunkl operator is the derivative operator and the Dunkl factorial becomes the classical factorial. 

In~\cite{GMV} the Sheffer-Dunkl sequences have been introduced via Umbral Calculus (see~\cite{Rom}). For this purpose, it is necessary the function
\begin{equation}\label{eq:Dunkl-exp}
\Ea(x)=\sum_{p=0}^{\infty}\frac{x^p}{\cc_{p,\alpha}},
\end{equation}
that plays the role of the exponential function, and the Dunkl translation
\begin{equation}
\label{eq:Dunkl-tras}
\tau_yf(x)=\sum_{p=0}^{\infty}\frac{y^p}{\cc_{p,\alpha}}\La^pf(x),
\end{equation}
that generalizes the classical translation
\[
f(x+y)=\sum_{p=0}^{\infty}\frac{y^p}{p!}f^{(p)}(x).
\]
Particular cases of Sheffer-Dunkl polynomials are the Appell-Dunkl polynomials that have been studied in~\cite{CDPV, CMV, DPV, JMC, Ros} or the discrete Appell-Dunkl polynomials introduced in~\cite{ELMV, GLMV}. 

Following the same arguments as in~\cite{GMV} changing the sequence of the Dunkl factorials $\gamma=\{\cc_{p,\alpha}\}_{p=0}^{\infty}$ by the sequence $m=\{m(p)\}_{p=0}^{\infty}$, the Dunkl exponential by the function
\begin{equation}\label{eq:m-exponential}
E_m(x)=\sum_{p=0}^{\infty}\frac{x^p}{m(p)},
\end{equation}
if the sequence $m$ guarantees convergence in the whole complex domain, the power series in (\ref{eq:m-exponential}) determines an entire function. In practice, the sequence $m$ is assumed to be a strongly regular sequence, following~\cite{thilliez}. This function is the main object constructing the solutions of a linear system of moment differential equations, described in the next section.

\section{Linear systems of Dunkl differential equations}\label{sec3}

Let $\gamma=\{\cc_{p,\alpha}\}_{p=0}^{\infty}$ be with $\alpha>-1$. In this section we give the general solution to a system
\begin{equation}\label{eq:D-sistema}
\partial_{\gamma}y(z)=\La y(z)=A y(z),
\end{equation}
with $A\in\C^{n\times n}$ a constant matrix and $y(z)=(y_1(z),\ldots,y_n(z))^T$. 

The function~\eqref{eq:Dunkl-exp} determined by the sequence of moments $\gamma$ is the entire function known as the Dunkl exponential and can be also expressed as 
\[
\Ea(z) = \I_\alpha(z) + \frac{1}{2(\alpha+1)}G_{\alpha}(z),
\]
where
$$\I_\alpha(z) = 2^\alpha \Gamma(\alpha+1)\frac{J_\alpha(iz)}{(iz)^\alpha} = \sum_{p=0}^\infty \frac{z^{2p}}{\cc_{2p,\alpha}}$$
and
\[
G_\alpha(z)\;=\;z\,\I_{\alpha+1}(z)
\;=\;\sum_{p=0}^\infty \frac{z^{2p+1}}{\cc_{2p,\alpha+1}}
\;=\;2(\alpha+1)\sum_{p=0}^\infty \frac{z^{2p+1}}{\cc_{2p+1,\alpha}}.
\]

Here, $J_\alpha(z)$ is the Bessel function of order $\alpha$ (and hence, $\I_\alpha(z)$ is a small variation of the modified Bessel function of the first kind, $I_\alpha(z)$). The Dunkl exponential satisfies 
\[
\La(\Ea(\lambda x)) = \lambda \Ea(\lambda x),
\]
and when applying the Dunkl translation~\eqref{eq:Dunkl-tras} to the Dunkl exponential we have
\begin{equation}\label{eq:D-sum}
\tau_y\Ea(x)=\Ea(x)\Ea(y).
\end{equation}
Note that~\eqref{eq:D-sum} generalizes the formula $e^{x+y}=e^xe^y$.

%

In order to obtain the solutions of \eqref{eq:D-sistema}, we need to consider the following formal power series
$$
\Eah(\lambda z)=\sum_{p\ge h}{p\choose h}\frac{\lambda^{p-h}z^p}{\cc_{p,\alpha}}, \quad \lambda\in\C,\, h\ge 0,
$$
that it is also an entire function for every $h\ge 0$.


We also need the Jordan canonical form of the matrix $A$. Let us consider a matrix $A\in\C^{n\times n}$ with eigenvalues $\{\lambda_j\}_{j=1}^m$ of multiplicity $\{\l_j\}_{j=1}^m$. Let
\[
v_{j,k}, \quad 1\le j \le m, \quad 1\le k \le \l_j,
\]
be linearly independent vectors such that
\begin{itemize}
\item $\sum_{j=1}^n \l_j=n$,
\item $A v_{j,1}=\lambda_j v_{j,1}$, $j=1$, \dots, $m$.
\item $Av_{j,k}=v_{j,k-1}+\lambda_j v_{j,k}$,  $j=1$, \dots, $m$, $k=2$, \dots $\l_j$.
\end{itemize}
Then, using the results of~\cite{lastra} we can prove
\begin{theorem}
Let $A\in\C^{n\times n}$ be a matrix whose eigenvalues are $\{\lambda_j\}_{j=1}^m$ of multiplicity $\{\l_j\}_{j=1}^m$. The functions
\[
y_j(z)=v_{j,k}\Ea(\lambda_jz)+v_{j,k-1}E_{\alpha,1}(\lambda_jz)+\cdots+v_{j,1}E_{\alpha,k-1}(\lambda_jz),\quad k=1,\ldots, \l_j,\quad j=1,\ldots,m,
\]
form a fundamental system of solutions of \eqref{eq:D-sistema}.
\end{theorem}

In the previous result, it has been omitted that the set of solutions to (\ref{eq:D-sistema}) is a vector space of dimension $n$. We refer to~\cite{lastra} for further details.

In~\cite{lastra}, the fact that the sequence $m$ defining moment derivative is a strongly regular sequence admitting a nonzero proximate order guarantees asymptotic properties of the entire solutions to the linear system of moment differential equations at infinity. This is the case of  the sequence $\gamma$. This guarantees some asymptotics properties for the solutions of \eqref{eq:D-sistema}. We start recalling the concept of strongly regular sequence, following~\cite{thilliez}.

\begin{defin}\label{defi132}
Let $\mathbb{M}:=\{M_p\}_{p=0}^{\infty}$ be a sequence of positive real numbers such that $M_0=1$. $\mathbb{M}$ is a strongly regular sequnce if it satisfies the following properties:
\begin{itemize}
\item[(lc)] $\mathbb{M}$ is logarithmically convex, i.e., $M_{p}^2\le M_{p-1}M_{p+1}$ for every $p\ge1$.
\item[(mg)] $\mathbb{M}$ is of moderate growth, i.e., there exists $A_1>0$ with $M_{p+q}\le A_1^{p+q}M_pM_q$ for every pair of integers $p,q\ge0$.
\item[(snq)] $\mathbb{M}$ is non-quasianalytic, i.e., there exists $A_2>0$ such that
$$\sum_{q\ge p}\frac{M_q}{(q+1)M_{q+1}}\le A_2\frac{M_p}{M_{p+1}}$$
for every $p\ge0$.
\end{itemize}
\end{defin}



%

\begin{proposition}\label{prop:D-SR}
Let $\alpha>-1$. The sequence of Dunkl factorials $\{\cc_{p,\alpha}\}_{p=0}^{\infty}$ is a strongly regular sequence if $-1<\alpha\le 0$.
\end{proposition}
\begin{proof}
We start proving that $\{\cc_{p,\alpha}\}_{p=0}^{\infty}$ is (lc) for $-1<\alpha\le 0$. If $p=2k$, we have to prove that
\[
\cc_{2k,\alpha}^2\le \cc_{2k-1,\alpha}\cc_{2k+1,\alpha}.
\]
This is equivalent to $\alpha>-1$. When $p=2k+1$, the inequality
\[
\cc_{2k+1,\alpha}^2\le \cc_{2k,\alpha}\cc_{2k+2,\alpha},
\]
holds if $\alpha\le 0$. 

On the other hand, it is easy to prove that the sequences of the quotients $\{\frac{\gamma_{p,\alpha}}{\gamma_{p-1,\alpha}}\}_{p=0}^\infty$ and $\{\frac{ p!}{(p-1)!} \}_{p=0}^\infty$ are equivalent in the sense of Definition 2.6 of \cite{JG-S}, that is, there exist constants $C_1,\,C_2>0$ such that
\[
C_1 p\le \theta_{p,\alpha}\le C_2 p.
\]
 So, from Remark 2.8 of~\cite{JG-S}, as the sequence $\{p!\}_{p=0}^\infty$ is strongly regular, the sequence of Dunkl factorials is also strongly regular.

\end{proof} 
Another important concept is that of proximate order.
\begin{defin}
A nonzero proximate order $\rho(t)$ is a nonnegative continuously differentiable function defined in an interval of the form $(c,\infty)$ for some $c\in\R$ such that
$$\lim_{r\to\infty}\rho(r)=\rho,$$
for some $\rho\in\R$, and 
$$\lim_{r\to\infty}r\rho'(r)\ln(r)=0.$$
\end{defin}

Let $m=\{m(p)\}_{p=0}^{\infty}$ be a sequence of positive numbers. Let $M:[0,\infty)\to [0,+\infty)$ (see~\cite{mandel}) be the function defined by $M(0)=0$ and
	\begin{equation}\label{e377}
		M(t):=\sup_{p\ge 0}\log\left(\frac{t^p}{m(p)}\right),\quad t>0.
	\end{equation}
Let us consider
\[
d(t):=\frac{\log M(t)}{\log(t)}.
\]
\begin{proposition}\label{prop:ord-aprox}
For the sequence of Dunkl factorials $\gamma=\{\cc_{p,\alpha}\}_{p=0}^{\infty}$, the function $d(t)$ is a nonzero proximate order.
\end{proposition}
\begin{proof}
It is easy to prove that
\[
\lim_{p\to\infty}\log\left(\frac{\theta_{p,\alpha}}{\cc_{p,\alpha}^{1/p}}\right)=1.
\]
Then, from Theorem 3.14 of \cite{JG-S}, $d(t)$ is a nonzero proximate order.
\end{proof}

We also point out that the number $\omega(\gamma)$ involved in the theory of generalized summability is determined by the limit
\[
\omega (\gamma )=\lim_{p\rightarrow \infty }\frac{\log \left( p+\frac{2\alpha +1}{2}\left( 1-\left( -1\right) ^{p}\right) \right) }{\log \left(p\right) }=1
\]
taking into account Theorem 1.2 in~\cite{JG-S}, and its consequences. In addition to that, it is straightforward to check that the functions $M(t)$ defined in \eqref{e377} for the sequences $\{\gamma_{p,\alpha}\}_{p=0}^{\infty}$ and $\{p!\}_{p=0}^{\infty}$ are comparable. Indeed, by Stirling’s formula one verifies that the  sequences $\{\gamma_{p,\alpha}\}_{p\ge0}$ and $\{p!\}_{p\ge0}$ are equivalent in the sense of Definition 2.5 of \cite{JG-S}, that is there exist constants $C_1,\,C_2>0$ such that
\[
C_1^p p!\le \gamma_{p,\alpha}\le C_2^p p!.
\]


\smallskip

In view of Proposition~\ref{prop:D-SR} and Proposition~\ref{prop:ord-aprox}, and regarding the theory of generalized summability (see Section 5,~\cite{sanz}), if $\alpha\in(-1,0]$, then there exists a pair of kernel functions for generalized summability $(e,E)$ associated to the sequence $\gamma$ such that
	\[
	\cc_{n,\alpha}=\int_0^{\infty} t^{n-1}e(t)\,dt,
	\]
	and the function
	\[
	E(z)=\sum_{p=0}^{\infty}\frac{z^p}{\cc_{p,\alpha}}
	\]
	is an entire function and there exist $k_3,\,k_4>0$ with
	\[
	|E(z)|\le k_3\exp\left(M\left(\frac{|z|}{k_4}\right)\right),\quad z\in \C.
	\]
Note that in this case, the function $E(z)$ is our Dunkl exponential $\Ea(z)$.


Additionally, in the recent paper~\cite{dugiva}, the authors prove the sequence $\gamma=\{\gamma_{n,\alpha}\}_{n=0}^{\infty}$ is the solution to a Hamburger moment problem, when $-1<\alpha<-1/2$. More precisely, it holds that
$$\gamma_{n,\alpha}=\int_{-\infty}^{\infty}t^n\omega_{\alpha}(t)dt,\quad n\ge0,$$
where 
$$\omega_{\alpha}(t)=\frac{|t|^{\alpha+1}(K_{\alpha}(|t|)+\hbox{sgn}(t)K_{\alpha+1}(|t|))}{2^{\alpha+1}\Gamma(\alpha+1)},$$
with $K_\alpha$ being the modified Bessel function of the second kind. A natural splitting of the previous construction is also shown in~\cite{dugiva}, writting for all $n\ge0$ 
$$\gamma_{n,\alpha}=(-1)^n\gamma_{n,\alpha}^{-}+\gamma_{n,\alpha}^{+},$$
where 
$$\gamma_{n,\alpha}^{-}=\int_0^{\infty}t^n\omega_{\alpha,-}(t)dt,\qquad \gamma_{n,\alpha}^{+}=\int_0^{\infty}t^n\omega_{\alpha,+}(t)dt,$$
for some completely monotonic functions $\omega_{\alpha,-},\omega_{\alpha,+}:(0,\infty)\to\R$ defined by
$$\omega_{\alpha,+}(t)=\frac{1}{2^{\alpha+1}\Gamma(\alpha+1)}t^{\alpha+1}(K_{\alpha}(t)+K_{\alpha+1}(t)),\qquad \omega_{\alpha,-}(t)=\frac{1}{2^{\alpha+1}\Gamma(\alpha+1)}t^{\alpha+1}(K_{\alpha}(t)-K_{\alpha+1}(t)),$$
for all $t>0$. We recall that Bessel functions of the second kind are multivalued functions for $\alpha$ not being an integer, and with a branch point at the origin.


Moreover, the theory of generalized summability leads to the fact that there exists $\beta>0$ such that for all $\theta\in(0,\pi)$ and $R>0$ there exists $k_5>0$ such that $|\Ea(z)|\le k_5|z|^{-\beta}$ for all $z\in\{z\in\C:|\hbox{arg}(z)-\pi|<\theta\pi\}\setminus D(0,R)$ where $D(0,R)$ is the open disc centered in zero and of radius $R$.
.

Being the sequence $\gamma$ a strongly regular sequence admitting a nonzero proximate order allows to apply the results in~\cite{lastra} on the asymptotic behavior of the solution to (\ref{eq:D-sistema}) near infinity. 

\begin{defin}
Let $f$ be an entire function in $\C$. We put $M_f(r)=\max\{|f(z)|:|z|=r\}$ for any $r\ge0$. The order of $f$ is defined by
$$\rho=\rho_f=\lim\sup_{r\to\infty}\frac{\ln^{+}(\ln^{+}(M_f(r)))}{\ln(r)},$$
with $\ln^+(\cdot)=\max\{0,\ln(\cdot)\}$. Given $f$ as above of order $\rho\in\R$, the type of $f$ is defined by
$$\sigma=\sigma_f=\lim\sup_{r\to\infty}\frac{\ln^{+}(M_f(r))}{r^\rho}.$$
\end{defin}

More precisely, one has the following concrete asymptotic results on the solution to (\ref{eq:D-sistema}).

\begin{proposition}[Theorem 5,~\cite{lastra}]
For $-1< \alpha\le 0$, any solution $y=y(z)$ of (\ref{eq:D-sistema}) satisfies that $y$ is a vector of entire functions of order 1, and with type upper bounded by $\sigma=\max\{|\lambda|:\lambda\in\hbox{spec}(A)\}$, or an entire function of order 0.
\end{proposition}

\begin{proposition}[Proposition 4,~\cite{lastra}]
For $-1< \alpha\le 0$, given a diagonalizable matrix $A\in\C^{n\times n}$, then any component $y_j$ of the solution $y=y(z)$ of (\ref{eq:D-sistema})  satisfies that
$$|y_j(re^{i\theta})|\le \frac{C}{r^{\beta}},\quad r\ge R_0,\quad 1\le j\le n,$$
for some $C,\beta>0$ and some $R_0>0$, provided that $\theta\in\mathcal{A}$, if $\mathcal{A}$ is nonempty, where
$$\mathcal{A}=\bigcap_{\lambda\in\hbox{spec}(A)}\left\{\theta\in\R:\frac{\pi}{2}-\hbox{arg}(\lambda)<\theta<\frac{3\pi}{2}-\hbox{arg}(\lambda)\right\}.$$

\end{proposition}

Following~\cite{maergoiz,maergoiz2}, given a nonzero proximate order $\rho=\rho(t)$ of an entire function $f$, the generalized indicator of $f$ is defined by
$$h_f(\theta)=\lim\sup_{r\to\infty}\frac{\ln|f(re^{i\theta})|}{r^{\rho(r)}},\quad \theta\in\R.$$

\begin{proposition}
For $-1<\alpha\le 0$, given a diagonalizable matrix $A\in\C^{n\times n}$, then any component $y_j$ of the solution $y=y(z)$ of (\ref{eq:D-sistema}) satisfies that
$$h_{y_j}(\theta)\le\max\{|\lambda|h_{E}(\theta+\hbox{arg}(\lambda)):\lambda\in\hbox{spec}(A)\},$$
for all $\theta\in\R$. Here, $h_{E}(\theta)=\cos(\theta)$ if $|\hbox{arg}(\theta)|\le\frac{\pi}{2}$ and $h_E(\theta)=0$ otherwise, for $\theta\in[-\pi,\pi]$.
\end{proposition}

\section{$m$-translation and application}\label{sec4}

Having shown that Dunkl operator can be seen as a particular realization of moment derivation when considering sequence of Dunkl factorials as the sequence of moments, we now lean on this theory to make advances in the more general, but also rather more unexplored, theory of moment functional equations.

A natural extension to Dunkl translation can be stated when considering a more general sequence of positive real numbers $m$. In the context of Dunkl operators, the Dunkl translation operator is defined  by (see~\cite{GMV})
$$\tau_yf(x)=\sum_{p=0}^{\infty}\Lambda_\alpha^pf(x)\frac{y^p}{\gamma_{p,\alpha}},\quad \alpha>-1.$$
When regarding more general sequences, one can generalize the previous definition as follows, initially from a formal point of view.

\begin{defin}
Let $y\in\C$ and let $m=\{m(p)\}_{p=0}^{\infty}$ be a sequence of positive real numbers. We define the operator $\tau_{y,m}:\C[[z]]\to\C[[z]]$ by 
\begin{equation}\label{eq:m-tras}
\tau_{y,m}\hat{f}(z)=\sum_{p=0}^{\infty}\frac{y^p}{m(p)}\partial_m^p\hat{f}(z).
\end{equation}
\end{defin}

Observe that $\tau_{y,m}$ recovers Dunkl translation $\tau_y$ when departing from the strongly regular sequence $m=\gamma$ of Dunkl factorials.

The previous definition can naturally be extended to holomorphic functions defined on some neighborhood of the origin via Taylor representation of the function at the origin.
We will denote by $\mathcal{O}(U)$ the set of holomorphic functions on the open $U$ and $D(0,R)$ will be the open disc centered in zero and of radius $R$.

\begin{proposition}
Let $m=\{m(p)\}_{p=0}^{\infty}$ be a strongly regular sequence and $f\in\mathcal{O}(D(0,R))$, for some $R>0$. For every $y\in D(0,R/A_1)$, one has that $\tau_{y,m}f$ defines a holomorphic function in $D(0,R/A_1)$, where $A_1$ is the constant appearing in (mg) condition. 
\end{proposition}
\begin{proof}
Let $f\in\mathcal{O}(D(0,R))$. Then, from Taylor expansion of $f$ at the origin, one has that 
$$f(z)=\sum_{p=0}^{\infty}\frac{\partial_m^pf(0)}{m(p)}z^p,$$
for all $z\in D(0,R)$.  

After applying Cauchy integral formula for derivatives, for any $0<R'<R$  we arrive at
$$|\partial_m^\ell f(0)|=\frac{m(\ell)}{\ell!m(0)}|f^{(\ell)}(0)|=\frac{m(\ell)}{2\pi  m(0)}\left|\oint_{|\omega|=R'}\frac{f(\omega)}{\omega^{\ell+1}}d\omega\right|,$$
valid for all $\ell\ge0$ and all $z\in D(0,R')$. Usual estimates yield the existence of $M>0$ such that
\begin{equation}\label{eq:acotaciones}
	|\partial_m^\ell f(0)|\le M\left(\frac{1}{R'}\right)^{\ell}m(\ell),\quad z\in D(0,R').
\end{equation}

Now, we are going to see that $\partial_m^nf$ defines an holomorphic function in $D(0,R/A_1)$ for $n\ge 0$.
More precisely, given $0<R'<R$, we choose $R'<R''<R$. Then, using \eqref{eq:acotaciones} and (mg) of Definition~\ref{defi132}, there exists $M>0$ such that for all $z\in D(0,R'/A_1)$ one has that

$$ |\partial_m^nf(z)|\le \sum_{p=0}^{\infty}\frac{|\partial_m^{p+n}f(0)|}{m(p)}|z|^{p}\le M\frac{1}{(R'')^n}\sum_{p=0}^{\infty}\frac{m(p+n)}{m(p)}\left(\frac{|z|}{R''}\right)^p\le Mm(n)\left(\frac{A_1}{R''}\right)^n\sum_{p=0}^{\infty}\left(\frac{|z|A_1}{R''}\right)^p,$$
which is convergent for $z\in D(0,R'/A_1)$.

Finally, we have for all $y\in\C$ with $|y|<R'/A_1$ and $z\in D(0,R'/A_1)$ that
$$|\tau_{y,m}f(z)|\le \sum_{n=0}^{\infty}\frac{|y|^n}{m(n)}\tilde{M}\frac{m(n)A_1^n}{(R'')^n}=\tilde{M}\sum_{n=0}^{\infty}\left(\frac{|y|A_1}{R''}\right)^n$$
with $\tilde{M}=M \sum_{p=0}^{\infty}\left(\frac{R'}{R''}\right)^p$. 
\end{proof}

Now, we introduce some necessary notation. Let $\mathcal{R}$ be the Riemann surface of the logarithm. Let $\theta>0$ and $d\in\R$. Then
\[
S_d(\theta)=\{z\in\mathcal{R}:\, |\arg(z)-d|<\theta/2\},
\]
is the open infinite sector contained in the Riemann surface of the logarithm with the vertex at the origin, bisecting direction $d\in\R$ and opening $\theta>0$. We write $S_d$ when the opening $\theta>0$ is not specified.  

Given an unbounded open sector $S_d\subseteq\C$ with vertex at the origin, a neighborhood of the origin $D$, and a sequence of positive real numbers $\mathbb{M}$. We denote by $\mathcal{O}^{\mathbb{M}}(S_d\cup D)$ the set of all holomorphic functions in $S_d\cup D$ such that for all infinite sector with vertex at the origin $V$ with $\overline{V}\setminus\{0\}\subseteq S_d$ and any neighborhood of the origin $B$ with $\overline{B}\subseteq D$ one has that there exist $c,\,k>0$ with
$$
|f(z)|\le c\exp\left(M\left(\frac{|z|}{k}\right)\right),\quad z\in V\cup B,
$$
where $M(t)$ is the function \eqref{e377} for the sequence $\mathbb{M}$.

Regarding the growth at infinity of the $m$-translation of a function with certain growth at infinity, one has the following result.

\begin{proposition}
Let $m=\{m(p)\}_{p=0}^{\infty}$ and $\tilde{m}$ be strongly regular sequences. Let $S_d$ be an infinite sector with vertex at the origin and bisecting direction $d\in\R$, and let $r>0$. Assume that $f\in \mathcal{O}^{\tilde{m}}(S_d\cup D(0,r))$. Then, there exist $R>0$ and $0<\tilde{r}<r$ such that for all infinite sector $V$ with $\overline{V}\setminus\{0\}\subseteq S_d$, the function $z\mapsto \tau_{y,m}f(z)$ belongs to $\mathcal{O}^{\tilde{m}}(V\cup D(0,r'))$ for all $y\in D(0,R)$.
\end{proposition}
\begin{proof}
In view of the previous proposition, one has that $\tau_{y,m}f\in\mathcal{O}(D(0,\tilde{r}))$ for all $y\in D(0,R)$, for some $R>0$. In view of (2) in Theorem 3,~\cite{lamisu}, one derives that for every infinite sector $V$ with vertex at the origin and $\overline{V}\setminus\{0\}\subseteq S_d$, one has that $\partial_m^{n}f\in\mathcal{O}(V\cup D(0,\tilde{r}))$. In addition to this, for all $z\in V\cup D(0,\tilde{r})$ one has that
$$|\tau_{y,m}f(z)|\le \sum_{n=0}^{\infty}\frac{|y|^n}{m(n)}|\partial^{n}_mf(z)|\le C_1\exp(\tilde{M}(C_3|z|)\sum_{n=0}^{\infty}(C_2|y|)^n$$
for some $C_1,\,C_2,\,C_3>0$, where $\tilde{M}$ is the function  \eqref{e377} associated to the sequence $\tilde{m}$. This yields the result, when choosing $y\in D(0,\min\{1/C_2,R/A_1\})$.
\end{proof}

\textbf{Remark:} We observe that summability is not closed under generalized translation. Let us consider the series
$$\hat{f}(z)=\sum_{n=0}^{\infty}(-1)^{n}n!z^n,$$
which is a classical example of 1-summable series along any direction $d\neq\pi$. However, $\tau_{y,m}\hat{f}$ is not $\mathbb{M}$-summable along any direction for any strongly regular sequence $\mathbb{M}=\{M_p\}_{p=0}^{\infty}$. Indeed, Nevanlinna's theorem states that  $\mathbb{M}$-summability along a direction of $\tau_{y,m}\hat{f}$ is equivalent to the fact that $\hat{\mathcal{B}}_{\mathbb{M}}\tau_{y,m}\hat{f}$ defines a holomorphic function on some neighborhood of the origin, which can be extended along an infinite sector with some precise growth at infinity. Here, $\hat{\mathcal{B}}_{\mathbb{M}}$ stands for $\mathbb{M}$ formal Borel operator sending a formal power series $\sum_{n=0}^{\infty}a_nz^n$ to $\sum_{n=0}^{\infty}\frac{a_n}{M_n}z^n$.

 Let us show that the formal power series $\hat{\mathcal{B}}_{\mathbb{M}}\tau_{y,m}\hat{f}$ has null radius of convergence for every $y\neq 0$. One has that
$$\hat{\mathcal{B}}_{\mathbb{M}}\tau_{y,m}\hat{f}=\sum_{n=0}^{\infty}(-1)^n\left(\sum_{p=0}^{\infty}\frac{y^p}{m(p)}(-1)^p(n+p)!m(n+p)\right)\frac{z^n}{M_nm(n)}.$$

Let $n\ge0$. Then, one has that the formal power series
$$\sum_{p=0}^{\infty}\frac{(n+p)!m(n+p)}{m(p)}(-y)^p$$
has null radius of convergence due to (lc) property satisfied by sequence $m$, which implies that $m(n+p)\ge m(n) m(p)$.

As a conclusion, we provide an application of the previous theory to solve certain families of functional equations involving $m$-translations, based on the results obtained by W. Strodt in~\cite{st}.

\begin{theorem}
Let $m$ be a strongly regular sequence and let $\omega_{\ell},\,c_{\ell}$ be complex numbers for $1\le \ell\le n$ for some positive integer $n$, with $c_j\neq 0$ for some $1\le j\le n$. We consider the functional equation 
\begin{equation}\label{e509}
c_1\tau_{\omega_1,m}y(z)+c_2\tau_{\omega_2,m}y(z)+\cdots+ c_n\tau_{\omega_n,m}y(z)=0.
\end{equation}
Then, the following statements hold:
\begin{itemize}
\item The set of entire solutions to (\ref{e509}) is a vector space over $\C$.
\item Let $E_m$ be the entire function defined by (\ref{eq:m-exponential}) and let $f(z)$ be the entire function given by
\begin{equation}\label{e530}
f(z)=c_1E_m(\omega_1 z)+\cdots+ c_nE_m(\omega_n z).
\end{equation}
For all $z_0\in \C$ such that $f(z_0)=0$, 
one has that the function 
$$y(z)=E_m(z_0z)$$
is an entire solution to (\ref{e509}).
\item Let $\{z_{0,1},\ldots,z_{0,N}\}$ be a set of $N$ different roots of $f$. Then, the set $\{E_m(z_{0,1}z),\ldots,E_m(z_{0,N}z)\}$ conforms a set of $N$ linearly independent solutions to (\ref{e509}). 
\end{itemize}
\end{theorem}
\begin{proof}
The first part is a consequence of the fact that the null function is an entire solution of (\ref{e509}) and  
$$\tau_{\omega,m}(c_1y_1+c_2y_2)= c_1\tau_{\omega,m}y_1+c_2\tau_{\omega,m}y_2,$$
for all entire functions $y_1,y_2$ and $c_1,c_2\in\C$. 

The second statement is a consequence of the following property.
\begin{lemma}
Let $\omega,\xi\in\C$ and let $m=\{m(p)\}_{p=0}^{\infty}$ be a strongly regular sequence. Then, it holds that
$$\tau_{\omega,m}E_m(\xi z)=E_m(\xi \omega) E_m(\xi z),\quad z\in \C.$$
\label{lemma-tauonE}
\end{lemma}
\begin{proof}
Direct inspection of $m$-translation operator yields
$$\tau_{\omega,m}E_m(\xi z)=\sum_{p=0}^{\infty}\frac{\omega^p}{m(p)}\partial_m^pE_m(\xi z).$$
We observe that for all $p\ge0$ one has that $\partial^p_m(E_m(\xi z))=\xi^pE_m(\xi z)$. Therefore,
$$\tau_{\omega,m}E_m(\xi z)= \sum_{p=0}^{\infty}\frac{\omega^p\xi^p}{m(p)}E_m(\xi z)=E_m(\xi \omega) E_m(\xi z).$$
\end{proof}

Let $z_0\in\C$ with $f(z_0)=0$. We observe that the function $z\mapsto E_m(z_0z)$ is an entire function. Moreover, it holds that
$$
c_1\tau_{\omega_1,m}E_m(z_0z)+\ldots+ c_n\tau_{\omega_n,m}E_m(z_0z)=\left(\sum_{j=1}^{n}c_jE_m(z_0 \omega_j)\right)E_m(z_0 z)=0,$$
for all $z\in\C$.

For the last statement, we fix a set $\{z_{0,1},\ldots,z_{0,N}\}$ of $N$ different complex numbers such that $f(z_{0,j})=0$ for $1\le j\le N$. Let us consider a linear combination of the functions $E_m(z_{0,j})$ which is null, i.e. 
$$\lambda_1E_m(z_{0,1}z)+\ldots+\lambda_NE_m(z_{0,N}z)=0.$$
Applying $\partial_m^{j}$, $0\le j\le N-1$, on the previous equality and evaluating at $z=0$ one arrives at a linear system of $N$ equations of the form
$$z_{0,1}^j\lambda_1+\cdots+z_{0,N}^j\lambda_N=0,\qquad 0\le j\le N-1.$$
The previous linear system has a Van der Monde matrix as its matrix of coefficients, with determinant given by $\prod_{1\le i<j\le N}(z_{0,i}-z_{0,j})\neq 0$. Therefore, the only solution to the previous system is given by $\lambda_1=\ldots=\lambda_N=0$, and the result follows.
\end{proof}

\begin{remark}
In view of the previous result, the problem of finding solutions to (\ref{e509}) is reduced to the problem of finding the zeros of an entire function of the form (\ref{e530}), which is a hard problem. Indeed, in the classical framework where $m$ is the sequence of factorial numbers the equation (\ref{e509}) becomes a classical difference equation
\[
c_1y(z+\omega_1)+c_2y(z+\omega_2)+\cdots+c_ny(z+\omega_n)=0,
\]
and $y(z)=e^{z_0z}$ is an entire solution where $z_0$ is a zero of the exponential polynomial
\[
f(z)=c_1e^{\omega_1z}+c_2e^{\omega_2z}+\cdots+c_ne^{\omega_nz}.
\]
The zeros of this kind of functions have been widely studied in literature. See~\cite{abha,ritt,sepvid,mor} among others.

When we take the sequence $\gamma$ of the Dunkl factorials, one solution of \eqref{e509} is given by
\[
y(z)=\Ea(z_0z),
\]
where $z_0$ is a zero of
\[
f(z)=c_1\Ea(\omega_1z)+c_2\Ea(\omega_2z)+\cdots+c_n\Ea(\omega_nz).
\]
 \end{remark}

\section{The $m$-even translation operator}



In this section we complement the Dunkl translation operator (\ref{eq:Dunkl-tras}) by revisiting its $m$-even counterpart. This operator was already introduced in \cite{GMV}, and we return to it here to place it within the general moment–derivative framework developed in this work and to record several consequences that will be needed later. The $m$-even translation is the natural companion to the $m$-translation, i.e., it is precisely the symmetrization of $\tau_{y,m}$ and therefore inherits, without additional hypotheses, all the properties established for $\tau_{y,m}$. Conceptually, it isolates the even component of shifts and furnishes a parity-preserving operator tailored to our framework.

\smallskip

Classically, for $f:\mathbb{R}\rightarrow \mathbb{C}$ and a step $y\in 
\mathbb{R}$, the even translation operator on $f$ is defined by 
\[
(T_{y}\,f)(x):=\frac{f(x+y)+f(x-y)}{2},\qquad x\in \mathbb{R}. 
\]%
Equivalently, $T_{y}=\tfrac{1}{2}(\tau _{y}^{\mathrm{cl}}+\tau _{-y}^{%
	\mathrm{cl}})$, where $\tau _{y}^{\mathrm{cl}}f(x)=f(x+y)$ is the classical
translation operator. Thus $T_{y}$ picks the even part of the shift by $\pm y
$, and admits the Taylor representation 
\[
(T_{y}\,f)(x)=\sum_{k=0}^{\infty }\frac{y^{2k}}{(2k)!}\,f^{(2k)}(x), 
\]%
so only the even derivatives of $f(x)$ contribute.

\smallskip

At this point it is worth noting that, beyond the classical merits of the
classical even translation operator, its $q$-analogue known as $q$-even
translation operator and related constructions arise naturally in quantum
calculus: they underlie Calderon-type reproducing formulas and the design of 
$q$-wavelets, where positivity and boundedness are key features (see, e.g., 
\cite{FBJNMP06,NSIM13} and the references therein). 

\smallskip 

Let $m=\{m(p)\}_{p\geq 0}$ be a strongly regular sequence as defined in \ref%
{defi132}, and denote by $\partial _{m}$ the associated moment derivative,
by $E_{m}$ the entire function in (\ref{eq:m-exponential}), and by $\tau
_{y,m}$ the $m$-translation operator in (\ref{eq:m-tras}). For $y\in \mathbb{%
	C}$, we define the $m$-even translation operator $T_{y,m}:\mathbb{C}%
[[z]]\rightarrow \mathbb{C}[[z]]$ by 
\[
\left( T_{y,m}\,f\right) (z):=\frac{\tau _{y,m}\,f(z)+\tau _{-y,m}\,f(z)}{2}.
\]

Under this framework, from (\ref{eq:m-tras}) one has 
\[
\tau _{y,m}\,f(z)=\sum_{n=0}^{\infty }\frac{y^{n}}{m(n)}\partial_{m}^{n}%
\,f(z)\text{\quad and\quad }\tau _{-y,m}\,f(z)=\sum_{n=0}^{\infty }\frac{%
	(-y)^{n}}{m(n)}\partial _{m}^{n}\,f(z), 
\]
and it is trivial to check that averaging cancels all odd $n$ and keeps the
even ones, yielding the following even-derivative expansion 
\[
(T_{y,m}\,f)(z)=\sum_{k=0}^{\infty }\frac{y^{2k}}{m(2k)}%
\partial_{m}^{2k}f(z), 
\]%
which holds for all functions $f(z)$, formal or holomorphic in a
neighborhood of the origin where $\tau_{y,m}$ is defined.

\smallskip

Concerning the spectral action on $E_{m}$, set $E_{m}^{\xi }(z):=E_{m}(\xi z)
$. By Lemma \ref{lemma-tauonE}, for all $\xi ,y,z\in \C$, 
\[
\tau _{y,m}\,E_{m}^{\xi }(z)=E_{m}(\xi y)\,E_{m}^{\xi }(z),\qquad \tau
_{-y,m}\,E_{m}^{\xi }(z)=E_{m}(-\xi y)\,E_{m}^{\xi }(z).
\]%
Averaging gives 
\[
\left( T_{y,m}E_{m}^{\xi }\right) (z)=\frac{E_{m}(\xi y)+E_{m}(-\xi y)}{2}%
\,E_{m}^{\xi }(z),
\]%
so $E_{m}^{\xi }$ diagonalizes $T_{y,m}$ with eigenvalue equal to the even
part of $E_{m}$ at $\xi y$.

\smallskip 

With respect to the Dunkl particularization, let $\alpha >-1$ and let $%
\gamma =\{\gamma _{p,\alpha }\}_{p\geq 0}$ be the Dunkl factorials as in (%
\ref{eq:ccna}), with Dunkl operator $\Lambda _{\alpha }$ and Dunkl
translation 
\[
\tau _{y,\alpha }\,f(z)=\sum_{n=0}^{\infty }\frac{y^{n}}{\cc _{n,\alpha }}%
\Lambda _{\alpha }^{n}\,f(z),\qquad y,z\in \mathbb{C}.
\]
The $m$-even translation introduced above specializes to the Dunkl even
translation (see \cite[Example 3.2]{GMV}) by%
\[
(T_{y,\alpha }\,f)(z):=\frac{\tau _{y,\alpha }\,f(z)+\tau _{-y,\alpha }\,f(z)%
}{2}=\sum_{k=0}^{\infty }\frac{y^{2k}}{\cc _{2k,\alpha }}\,\Lambda_{\alpha }^{2k}f(z),
\]%
where the averaging cancels all odd powers as in the general $m$-case. Let
us denote next the Dunkl exponential by 
\[
E_{\alpha }(w)=\sum_{p=0}^{\infty }\frac{w^{p}}{\cc _{p,\alpha }}
=\I_{\alpha }(w)+\frac{1}{2(\alpha +1)}G_{\alpha }(w),
\]
with $\I_{\alpha }$ even and $G_{\alpha }$ odd. Then $E_{\alpha}(w)+E_{\alpha }(-w)=2\,\I_{\alpha }(w)$. Using the multiplicative law $\tau
_{y,\alpha }E_{\alpha }(\xi z)=E_{\alpha }(\xi y)\,E_{\alpha }(\xi z)$, one
obtains the spectral action 
\[
\left( T_{y,\alpha }E_{\alpha }(\xi \cdot )\right) (z)=\frac{E_{\alpha }(\xi
	y)+E_{\alpha }(-\xi y)}{2}\,E_{\alpha }(\xi z)=\I_{\alpha }(\xi y)\,E_{\alpha
}(\xi z),
\]
so $E_{\alpha }(\xi \,\cdot )$ diagonalizes $T_{y,\alpha }$ with eigenvalue given by the even Dunkl--Bessel component $\I_{\alpha}(\xi y)$.

\section{Declarations}

\noindent\textbf{Conflict of interest.} The authors declare to have no conflict of interest to disclose.

\vspace{0.3cm}

\noindent\textbf{Data availability.} We have no data.

\vspace{0.3cm}

\noindent\textbf{Acknowledgements.}

The work of Judit M\'inguez Ceniceros has been partially
supported by grant PID2024-155593NB-C22 funded by MICIU/AEI of Spain.

The work of Alberto Lastra is partially supported by the project PID2022-139631NB-I00 of Ministerio de Ciencia e Innovaci\'on, Spain.

EJH acknowledges Spanish \textquotedblleft Agencia Estatal de Investigaci\'on\textquotedblright\ research project [PID2024-155133NB-I00], Ortogonalidad, Aproximaci\'{o}n e Integrabilidad: Aplicaciones en Procesos Estoc\'{a}sticos Cl\'{a}sicos y Cu\'{a}nticos.




\begin{thebibliography}{99}
\bibitem{abha} H. Abbas, A. Hajj-Diab, On the distribution of zeros of exponential polynomials and Shapiro conjecture. J. Appl. Funct. Anal. 6, No. 1 (2011) 48--56.
\bibitem{bayo} W. Balser, M. Yoshino, Gevrey order of formal power series solutions of inhomogeneous partial differential equations with constant coefficients. Funkcial. Ekvac. 53 (2010) 411--434.




\bibitem{CDPV}
\'O. Ciaurri, A. J. Dur\'an, M. P\'erez, J. L. Varona,
Bernoulli-Dunkl and Apostol-Euler-Dunkl polynomials with applications to series involving
zeros of Bessel functions.
J. Approx. Theory
235 (2018), 20--45.

\bibitem{CMV}
\'O. Ciaurri, J. M\'inguez Ceniceros, J. L. Varona,
{Bernoulli-Dunkl and Euler-Dunkl polynomials and their generalizations}.
Rev. R. Acad. Cienc. Exactas F\'is. Nat. Ser. A Mat. RACSAM.
113 (2019), 2853--2876.


\bibitem{Du}
C. F. Dunkl,
{Differential-difference operators associated to reflection groups}.
Trans. Amer. Math. Soc.
{311} (1989), 167--183.

\bibitem{DPV}
A. J. Dur\'an, M. P\'erez, J. L. Varona,
Fourier-Dunkl system of the second kind and Euler-Dunkl polynomials,
{J. Approx. Theory}
{245} (2019), 23--39.


\bibitem{dugiva} A. J. Dur\'an, A. Gil Asensi, J. L. Varona, {Moment problems associated to the Dunkl factorial}, Proc. Amer. Math. Soc., to appear 2025.

\bibitem{ELMV}
J. I. Extremiana Aldana, E. Labarga, J. M\'inguez Ceniceros, J. L. Varona,
{Discrete Appell-Dunkl sequences and Bernoulli-Dunkl polynomials of the second kind},
{J. Math. Anal. Appl.} 
{507} (2022), 125832, 20~pp. 

\bibitem{FBJNMP06} A. Fitouhi, N\'{e}ji Bettaibi, {Wavelet Transforms in Quantum Calculus}, Journal of Nonlinear Mathematical Physics {13}(4) (2006), 492--506. 

\bibitem{gs} I. M. Gelfand, G. E. Shilov, {Generalized Functions.} Vol. 2, Space of Fundamental and Generalized Functions, Academic Press, New York, 1965.

\bibitem{GLMV} A. Gil Asensi, E. Labarga, J. M\'inguez Ceniceros,  J. L. Varona, {Boole-Dunkl polynomials and generalizations},
{Rev. Real Acad. Cienc. Exactas F\'is. Nat. Ser. A-Mat. RACSAM} 
{118} (2024), Paper No.~16, 18~pp.



\bibitem{GMV}
A. Gil Asensi, J. M\'inguez Ceniceros, J. L. Varona,
{Sheffer-Dunkl Sequences Via Umbral Calculus in the Dunkl Context},
{Bull. Malays. Math. Sci. Soc.}
{47} (2024), no.6, Paper No. 172.

\bibitem{ichmich} K. Ichinobe, S. Michalik, {On the summability and convergence of formal solutions of linear $q$-difference-differential equations with constant coefficients}, Math. Ann. 389, No. 2 (2024) 1099--1130. 

\bibitem{JG-S} J. Jim\'enez-Garrido, J. Sanz, {Strongly regular sequences and proximate orders,} J. Math. Anal. Appl. 438,  (2016) 920--945.
 
\bibitem{jimisasc} J. Jim\'enez-Garrido, I. Miguel-Cantero, J.  Sanz, G. Schindl, {Optimal flat functions in Carleman-Roumieu ultraholomorphic classes in sectors}, Result. Math. 78, No. 3, Paper No. 98, 35 p. (2023). 
\bibitem{jisasc} J. Jim\'enez-Garrido, J. Sanz, G. Schindl, {Surjectivity of the asymptotic Borel map in Carleman-Roumieu ultraholomorphic classes defined by regular sequences}, Rev. R. Acad. Cienc. Exactas F\'is. Nat., Ser. A Mat., RACSAM 115, No. 4, Paper No. 181, 18 p. (2021). 
\bibitem{lasmich} A. Lastra, S. Michalik, {On sequences preserving $q$-Gevrey asymptotic expansions}, Anal. Math. Phys. 14, No. 2, Paper No. 17, 25 p. (2024). 
\bibitem{lamisu} A. Lastra, S. Michalik, M. Suwi\'nska, {Summability of formal solutions for some generalized moment partial differential equations.} Result. Math. 76, No. 1 (2021) Paper No. 22. 
\bibitem{lamisu2} A. Lastra, S. Michalik, M. Suwi\'nska, {Estimates of formal solutions for some generalized moment partial differential equations,} J. Math. Anal. Appl. 500 (2021), no. 1.
\bibitem{lamisu3} A. Lastra, S. Michalik, M. Suwi\'nska, {Summability of formal solutions for a family of generalized moment integro-differential equations,} Fract. Calc. Appl. Anal. 24 (2021) 1445--1476.
\bibitem{lamisu4} A. Lastra, S. Michalik, M. Suwi\'nska, {Multisummability of formal solutions for a family of generalized singularly perturbed moment differential equations,} Results Math. 78(2) (2023) 49.
\bibitem{laspri} A. Lastra, C. Prisuelos-Arribas, {Solutions of linear systems of moment differential equations via generalized matrix exponentials}, J. Differ. Equations 372 (2023) 591--611. 
\bibitem{lastra} A. Lastra, {Entire solutions of linear systems of moment differential equations and related asympotitc growth at infinity,} 
Differ. Equ. Dyn. Syst. 32 (2024), no. 4, 943--964.

\bibitem{maergoiz} L. S. Maergojz, {Indicator diagram and generalized Borel-Laplace transforms for entire functions of a given proximate order}, St. Petersbg. Math. J. 12, No. 2 (2001) 191--232; translation from Algebra Anal. 12, No. 2  (2000) 1--63. 
\bibitem{maergoiz2} L. S. Maergojz, Asymptotic characteristics of entire functions and their applications in mathematics and biophysics. Mathematics and its Applications (Dordrecht) 559. Dordrecht: Kluwer Academic Publishers, xxiv, 2003. 
\bibitem{mandel} S. Mandelbrojt, {S\'eries Adh\'erentes, Regularisation des suites, Applications}, Gauthier-Villars, Paris, 1952.
\bibitem{mi} S. Michalik, {Analytic solutions of moment partial differential equations with constant coefficients,} Funkcial. Ekvac. 56 (2013), no. 1, 19--50.
\bibitem{michalik12} S. Michalik, {Multisummability of formal solutions of inhomogeneous linear partial differential equations with constant coefficients}, J. Dyn. Control Syst. 18 (2012) 103--133.
\bibitem{misutk} S. Michalik, M. Suwi{\'n}ska, B. Tkacz, {On sequences preserving summability}, Result. Math. 80, No. 4, Paper No. 114, 30 p. (2025). 
\bibitem{mitk} S. Michalik, B. Tkacz, {The Stokes phenomenon for some moment partial differential equations,} J. Dyn. Control Syst. 25 (2019), no. 4, 573--598.


\bibitem{JMC}
J. M\'inguez Ceniceros,
Some Appell-Dunkl sequences,
{Bull. Malays. Math. Sci. Soc.} 
{46} (2023), Paper No.~64, 18~pp.

\bibitem{mor} C. J. Moreno, {The zeros of exponential polynomials. I}, Compos. Math. 26 (1973) 69--78.



\bibitem{NSIM13} A. Nemri, B. Selmi, {Calderón type formula in Quantum calculus}, Indagationes Mathematicae {24}(3) (2013), 491--504. 


\bibitem{remybook} P. Remy, {Asymptotic expansions and summability. Application to partial differential equations.} Lecture Notes in Mathematics 2351. Cham: Springer. xiii, 2024. 

\bibitem{Rom}
S. Roman,
{The Umbral Calculus},
Academic Press, Orlando, 1984.

\bibitem{Ros}
M. Rosenblum,
{Generalized Hermite polynomials and the Bose-like oscillator calculus},
{Oper. Theory Adv. Appl.}
{73} (1994), 369--396.

\bibitem{ritt} J. F. Ritt, {On the zeros of exponential polynomials}, Transactions A. M. S. 31 (1929) 680--686. 

\bibitem{sanz} J. Sanz, {Flat functions in Carleman ultraholomorphic classes via proximate orders}, J. Math. Anal. Appl. 415(2) (2014), 623--643.
\bibitem{sanzrev} J. Sanz, {Asymptotic analysis and summability of formal power series,} Analytic, algebraic and geometric aspects of differential equations, 199--262, Trends Math., Birkh\"auser/Springer, Cham, 2017.
\bibitem{sepvid} J. M. Sepulcre, T. Vidal, {Equivalence classes of exponential polynomials with the same set of zeros}, Complex Var. Elliptic Equ. 61, No. 2 (2016) 225--238. 
\bibitem{st} W. Strodt, {Linear difference equations and exponential polynomials}, Trans. Am. Math. Soc. 64 (1948) 439--466. 
\bibitem{su} M. Suwi\'nska, {Gevrey estimates of formal solutions for certain moment partial differential equations with variable coefficients,} J. Dyn. Control Syst. 27, No. 2 (2021) 355--370. 
\bibitem{thilliez} V. Thilliez, {Division by flat ultradifferentiable functions and sectorial extensions}, Result. Math. 44 (2003), 169--188.
\end{thebibliography}
\end{document}